\documentclass[11pt]{amsart}
\usepackage{amscd}
\usepackage{amsmath,amssymb,latexsym}
\usepackage{pb-diagram}
\usepackage{enumerate}
\usepackage{amsfonts}
\usepackage{amsthm}

\newtheorem{Thm}[equation]{Theorem}
\newtheorem{Cor}[equation]{Corrolary}

\newtheorem{Prop}[equation]{Proposition}
\newtheorem{Lem}[equation]{Lemma}
\newtheorem{Remark}[equation]{Remark}
\newtheorem{Def}[equation]{Definition}
\newcommand{\ra}{\rangle}
\newcommand{\la}{\langle}

\newcommand{\C}{\mathbb C}

\newcommand{\mR}{\mathcal R}

\newcommand{\bH}{\mathbb H}

\newcommand{\mS}{\mathcal S}
\newcommand{\mN}{\mathcal N}

\newcommand{\mF}{\mathcal F}

\newcommand{\mM}{\mathcal M}
\newcommand{\mO}{\mathcal O}

\newcommand{\mg}{\mathfrak g}

\newcommand\T{\operatorname{Tr}}
\newcommand\N{\operatorname{Nm}}

\newcommand\Span{\operatorname{Span}}

\newcommand\ind{\operatorname{ind}}
\newcommand\Ind{\operatorname{Ind}}

\newcommand\Hom{\operatorname{Hom}}

\theoremstyle{plain}
\title[Fourier transform on a cone]{Fourier transform on a cone and the minimal representation
  of even orthogonal group}
\author{Nadya Gurevich and David Kazhdan}
\address{School of Mathematics,
Ben Gurion University of the Negev, POB 653, Be'er Sheva 84105, Israel}
\address{Einstein Institute of Mathematics,
  The Hebrew University of Jerusalem, Givat Ram, Jerusalem, 9190401, Israel}
\email{ngur@math.bgu.ac.il}
\email{kazhdan@math.huji.ac.il}
\numberwithin{equation}{section}

\begin{document}
\maketitle
\begin{abstract} Let $G$ be an even orthogonal
quasi-split group defined over a local non-archimedean field $F$.
We describe the subspace of smooth vectors of the minimal representation
of $G(F),$ realized on the space of square-integrable functions on a  cone.
Our main tool is the Fourier transform on the cone, for which
we give an explicit formula. 
\end{abstract}  
\section{Introduction}\label{sec:intro}

\label{intro:part1}
\subsection{Notation}
\begin{itemize}
\item Let $F$ be  a local non-archimedean
  field, with the norm $|\cdot|$,
  the ring of integers $\mO_F$ and a fixed uniformizer
  $\varpi$ such that $|\varpi|=q^{-1},$
  where $q$ is the cardinality of the residue field.
  We fix an additive complex character $\psi$ of $F$.
  
\item Throughout this paper we use boldface characters for group
  schemes over $F$, such as ${\bf H}$, and plain text characters
  for their group of $F$-points, such as $H$. 
 \item 
  Let $({\bf V},q)$ be a non-degenerate quadratic space of dimension $2n+2$,
  $n\ge 3 $ and  the Witt index at least $n$, over $F$. We denote
  by $\la\cdot, \cdot \ra$ the associated bilinear form.
  The discriminant $\kappa\in F^\times/(F^\times)^2$ of $V$ gives rise to a
  quadratic algebra $K$ over $F$.
  The quadratic character $\chi_K$ is associated to this algebra
via class field theory. We also write $\chi_0$ for the trivial character of $F^\times$.
\item
We fix the decomposition $V=\bH\oplus V_1,$ where $\bH$ is a hyperbolic
plane and $V_1$ is a non-degenerate quadratic space of dimension $2n$.
\item
The group ${\bf G}={\bf O}({\bf V},q)$ is quasi-split. It
contains a maximal parabolic subgroup ${\bf Q}={\bf M}\cdot {\bf N},$
where ${\bf M}$ is canonically isomorphic to
${\bf GL_1\times O(V_1)}.$
Any  character of $GL_1$  is extended trivially to $M$.
\item
 We denote by
 $\la \cdot, \cdot\ra: N\times \bar N\rightarrow F$
 the non-degenerate pairing  given by the Killing form. 
 Fix a basis $\{e,e^\ast\}$ in $\bH$ of isotropic vectors
 such that $\la e, e^\ast\ra=1$ and  $N$ fixes the vector $e^\ast$.
 This gives rise to an isomorphism
 $$\varphi:\bar N\simeq V_1, \quad
 \varphi(\bar n)=proj_{V_1}(\bar n e^\ast),$$
 where $\bar N$ is the radical of the opposite parabolic $\bar Q$.
\item
 The characters of $N$ stay in
bijective correspondence with the vectors in $\bar N\simeq V_1$ by
$$\bar n \mapsto \Psi_{\bar n}(n)=\psi(\la n,\bar n\ra).$$
We shall call the character $\Psi_{\bar n}$ anisotropic
if the vector $\bar n\in V_1$ is anisotropic. 
\item 
  Denote by $C\subset \bar N$
  the cone  of isotropic vectors in $V_1\simeq \bar N$.
  The group $G_1=O(V_1)$ acts transitively on  the set
  $C_0=C\backslash \{0\}$. The action is on the right.
 We define an action of $Q$ on the space
  $\mS^\infty(C_0)$  of smooth functions on $C_0$  by:
\begin{equation}\label{intro:Q:action}
\left\{\begin{array}{ll}
%m\cdot f(w)=\chi_K(m)|m|^{n-1}f(m^{-1}wm) & m\in M\\
(a,h)\cdot f(w)=\chi_K(a)|a|^{n-1}f(awh) & a\in GL_1, h\in O(V_1),\\
n\cdot f(w)=\psi(-\la n, w\ra) f(w) & n\in N
\end{array}\right.
\end{equation}
\end{itemize}

%$a\cdot f(w)=\chi_K(a)|a|^{n-1}f(aw), \quad a\in GL_1$

\subsection{Minimal representations}
We show in section \ref{sec:minimal}  that the normalized principal
series $\Ind^{G}_{Q} \chi_K|\cdot|^{-1}$ contains unique
irreducible representation, denoted by $\Pi$. We prove
that $\Pi$ is minimal. The same argument applied to $G_1=O(V_1,q)$
gives the minimal representation $\Pi_1$ of $G_1$.

The realization of $\Pi$ on a space of functions
on  $C_0$ has been constructed by Savin \cite{Savin}, see also
\cite{SavinWoodbury}, for the case $G$ is split.
The result is extended to the case  of quasi-split $G$ 
with essentially the same proof, which we write in subsection
\ref{subsec:schrodinger}.

 \begin{Thm}\label{savin:intro}
  \begin{enumerate}
  \item
    There exists a $Q$ equivariant embedding
    $\Pi\hookrightarrow  \mS^\infty(C_0).$ We denote its image
    by $\mS$.
  \item  The space $\mS$ contains the space
    $\mS_c(C_0)$ of smooth functions of compact support
    and is contained in the space of smooth functions of bounded
    support.
\item
    The  Jacquet module $j_N(\mS)$ is isomorphic as $M$-module
    to the quotient $\mS/\mS_c(C_0)$, or equivalently to the
    space of germs  $\{[f]_0, f\in \mS\}$ at zero.  
  \end{enumerate}
  \end{Thm}

We seek to describe explicitly the representation of the group $M$ on the space 
$j_N(\mS)$. 
\iffalse
Let $\mS(|\cdot|^{-1})\subset \mS^\infty(C_0) $ be the 
subspace of functions satisfying
\begin{equation}
  \label{intro:ind}
f(aw)=\chi_K(a)|a|^{2-n}f(w), \quad a\in F^\times, w\in C_0.
\end{equation}
\fi
%
Consider the subspace of functions in $\mS^\infty(C_0)$
satisfying ${a\cdot f=|a|f}$ for all $a\in GL_1$. It is naturally isomorphic
as $M$-module to $|\cdot|\otimes \Ind^{G_1}_{Q_1}\chi_K |\cdot|^{-1}$
and hence  contains unique irreducible $M$-submodule $\mS_1,$ which is 
isomorphic to $|\cdot|\otimes \Pi_1$.

The goal of this paper is to construct an explicit isomorphism
 of $M$-modules
\begin{equation}\label{def:B}
  B: j_N(\mS) \xrightarrow{\sim} \mS_1 \oplus \C,
  \end{equation}
 where $M=GL_1\times O(V_1)$ acts on $\C$ by $\chi_K\,|\cdot|^{n-1}\otimes 1$.
\vskip 5pt
 The existence of such isomorphism  for the case
 $G$ is split is proven in \cite{Savin} using computations in the
 Iwahori-Hecke algebra of $G$ that are not readily
 generalized to the quasi-split case.

 \subsection{Fourier transform}
We make use of the operator $\Pi(r)\in Aut(\mS)$, where  
 $r\in O(V)$ is the involutive  element such that
\begin{equation}\label{def:r}
  r(e)=e^\ast, \quad r(e^\ast)=e, \quad r|_{V_1}=Id.
  \end{equation}
The operator $\Pi(r)$ is a unitary involution commuting with
$O(V_1)$ and is called for that reason {\sl a Fourier transform on the cone.}

In Theorem \ref{Pi(r)=Phi} we obtain an explicit formula
for the restriction of the operator $\Pi(r)$ to $\mS_c(C_0).$
Precisely, for $f\in \mS_c(C_0)$ one has $\Pi(r)(f)=\Phi(f),$ where 
\begin{equation}
  \Phi(f)( w)=
\gamma(\chi_K,\psi)
\int\limits_{F^\times}\hat \mR(f)(x w)
\psi(x^{-1}) \chi_K(-x) |x|^{n-2} d^\times x,
 \end{equation}
where $\gamma(\chi_K,\psi)$ is the Weil factor and
$\hat \mR$ is the normalized Radon transform,
(see section \ref{sec:radon}). 

To construct the isomorphism $B$ we proceed in the following steps:
\begin{itemize}
  \item
    Using the formula for $\Phi$ we write  the operator
$\Pi(r)|_{\mS_c(C_0)}=\Phi$ as a sum $\Phi_1 + \Phi_2$, such that
(see Proposition \ref{images:Phi12})  the image of $\Phi_1$
equals $\mS_1$. Since $\mS_1$ consists of homogeneous functions, 
the function $\Phi_1(f)$ is completely determined by its germ $[\Phi_1(f)]_0$
at zero. Besides,  the image of $\Phi_2$  is contained in the
space of locally constant functions on $C$. In particular,
$[\Phi_2(f)]_0$ is constant for any $f\in \mS_c(C_0)$. 
\item
We show  in Lemma \ref{S=S0+Phi(S0)} that
$\mS=\mS_c(C_0)+\Pi(r)(\mS_c(C_0)).$ 
Hence any vector in $j_N(\mS)$ is represented by
a germ of the form $[\Pi(r)f]_0$ for some  $f\in \mS_c(C_0)$.
\item All the above implies that the map
  $$B: j_N(\mS)\rightarrow \mS_1\oplus \C, \quad B([\Pi(r)f]_0)=
  (\Phi_1(f), [\Phi_2(f)]_0)$$
  is the isomorphism of $M$-modules.
\end{itemize}

Let us outline the content of the paper. The section $2$ is devoted
to the definition and properties of the minimal representation of $G$, $G_1$
and the model of  $\Pi$ in a space $\mS$ of smooth functions on the cone. 
The section $3$ is devoted to  Radon transform on the cone,
its equivariance properties and  asymptotics. The Fourier transform on the cone
is defined in terms of Radon transform. Finally, in section $4$ we
show that $\Pi(r)=\Phi$ on the space $\mS_c(C_0)$.

\begin{Remark}
  The space $(\Pi,\mS)$ is the subspace of smooth vectors
in the unitary  irreducible minimal representation
$(\hat \Pi, G, L^2(C)).$
 The operator $\hat\Pi(r)$ on $L^2(C)$
 is a particular example of the normalized intertwining  operator
 that was constructed in a general setting by Braverman and Kazhdan
 for split groups  \cite{BravermanKazhdan}. In particular,
 the space $\mS$ of the minimal
 representation coincides with the Schwartz space defined in {\sl loc.cit.}
\end{Remark}

\begin{Remark}\begin{enumerate}
    \item
  The Fourier transform acting on  the space $L^2(C)$, 
  for $F=\mathbb R$ has been studied by Kobayashi  and Mano
  \cite{KobayshiMano}.
\item For the case $n=3$ the space $C_0$ can be identified
  with the basic affine space of quasi-split group $SU_3$
  of $F$-rank $1$. In a forthcoming paper \cite{GurevichKazhdanBAS}
  we use the operator
  $\Phi$ to construct a family of generalized Fourier transforms
  on the basic affine space of any quasi-split simply-connected group.
\end{enumerate}
  \end{Remark}

\subsection*{Acknowledgment}
The research of the second author is
partially supported by the ERC grant No $669655$.
We thank the anonymous referee for his/her careful reading and questions
which helped to improve presentation and to remove inaccuracies
in the first version. 

\section{The degenerate principal series and the minimal representation}
\label{sec:minimal}

Let $\Pi_1$ be  a smooth minimal representation of $G_1=O(V_1)$.
This means that its character $\Theta_{\Pi_1}$,
viewed as an invariant distribution
on $G_1$, has an asymptotic Harish-Chandra-Howe expansion in a
neighborhood of $0$ in the Lie algebra $\mg_1$ of $G_1$:
$$\Theta_{\Pi_1}=c_1\hat\mu_{\mO_{min}}+c_0,$$
where $\mO_{min}$ is the minimal coadjoint orbit. 
In particular, its wave front $WF(\Pi_1)$ consists of the minimal orbit.
See \cite{HarishChandra} for details. 

Using the definition it is not easy to construct such representation
or even to know whether it exists. The idea is to look for $\Pi_1$
as a constituent of a well-understood representation, whose wave front
is already small. For orthogonal group
such representation is a degenerate principal series induced from a parabolic
subgroup $Q_1$, whose unipotent radical $N_1$ is abelian. 

\subsection{The degenerate principal series on $G_1=O(V_1)$}\label{subsec:dps}
We write the orthogonal decomposition $V_1=\bH\oplus V_2$ and fix a basis
$\{e_1,e_1^\ast\}$ of $\bH$  consisting of anisotropic vectors with $\la e_1,e_1^\ast\ra=1$. 
The parabolic subgroup $Q_1$ stabilizing the line $Fe_1^\ast$ admits a Levi
decomposition $M_1\cdot N_1$ with $M_1=GL_1\times O(V_2)$. 

Any complex character $\chi$ of $GL_1$ is canonically
extended to a character of $M_1$.

For a unitary character $\chi$ denote by $\chi_s$ the character
$\chi|\cdot|^s$ and by $I_{2n}(\chi,s)$
the normalized induced representation $\Ind^{G_1}_{Q_1}\chi_s.$
The standard intertwining operator
$\mM(\chi,s)\in \Hom_{G_1}(I_{2n}(\chi,s),I_{2n}(\chi^{-1},-s))$
corresponds to the shortest Weyl element $w$ such that
$w(Q_1)=\bar Q_1$, the opposite parabolic subgroup.
It is holomorphic for $\Re(s)\ge 0$.

\iffalse
We denote by $I_{2n}(\chi)$  the normalized induced representation
$\Ind^{G_1}_{Q_1}\chi$ and  by $\mM(\chi)$ the standard intertwining operator
corresponding to the shortest Weyl element $w$ such that $w(Q_1)=\bar Q_1$,
the opposite parabolic subgroup.
Thus $\mM(\chi)\in \Hom_{G_1}(I_{2n}(\chi),I_{2n}(\chi^{-1})).$
\fi

\begin{Prop} \label{red:ind}
  The representation $I_{2n}(\chi,s)$ is irreducible, unless
$\chi=\chi_0, s=\pm(n-1)$ or $\chi=\chi_K, s=\pm 1$.
\end{Prop}

  \iffalse
\item Let $\chi=|\cdot|^{1-n}$.
  The unique irreducible submodule of $I_{2n}(\chi)$
  is the trivial representation. The space
  $\Hom_{G_1}(I_{2n}(\chi^{-1}), I_{2n}(\chi))$ is one dimensional
 \item 
Let $\chi=\chi_K|\cdot|^{-1}$.
The  unique irreducible submodule $\Pi_1$  of $I_{2n}(\chi)$
is self-contragredient and unitarizable. 

The space $\Hom_{G_1}(I_{2n}(\chi^{-1}), I_{2n}(\chi))$
is one-dimensional and the image of any non-zero map is
the unique factor $\Pi_1$ of $I_{2n}(\chi)$.
\item
The  representation $\Pi_1$ is self-contragredient and unitarizable. 
\end{enumerate}
\end{Prop}
\fi
The reducibility of degenerate principal series of split, connected,
simply connected, simply-laced groups induced from maximal parabolic
subgroups with abelian radical has been studied by Weissman in \cite{Weissman}
using the Fourier-Jacobi functor. The method is also applicable
to the case $G_1=O(V_1)$ with minor changes. When the proof of
\cite{Weissman} goes through without changes we refer to it. 
Before proving Proposition \ref{red:ind} we define the Fourier-Jacobi
functor and record its useful properties.

Consider an orthogonal  decomposition $V_2=\bH\oplus V_3$,
and fix a basis $\{e_2,e_2^\ast\}$ of $\bH$ of isotropic vectors
such that $\la e_2, e_2^\ast\ra=1$. 
Put $X=Span\{e_1,e_2\}$ and $X^\ast=Span\{e_1^\ast,e_2^\ast\}$.
Both spaces $X,X^\ast$ have a  structure of symplectic space, with 
the symplectic form satisfying $[e_1,e_2]=[e_1^\ast,e_2^\ast]=1$. 

 The Heisenberg parabolic subgroup $P_1$ of $G_1,$ stabilizing the
 subspace $X^\ast$ in $V_1,$ admits Levi decomposition
  $P_1\simeq (GL(X^\ast)\times O(V_3))\cdot H$, where $H$ is a Heisenberg subgroup
  with center $Z$. 
  The natural map $H\rightarrow \Hom(V_3,X^\ast)$ has $Z$ as a  kernel.
  Further, $H/Z\simeq V_3^\ast\otimes X^\ast$ has a structure of a symplectic space. 

The natural homomorphism $SL(X^\ast)\times O(V_3)\rightarrow Sp(H/Z)$
splits in the double  covering of $Sp(H/Z)$, which gives rise to the  Weil representation
$\omega_{\psi,q}$ of  $(SL(X^\ast)\times O(V_3))\cdot H.$
We shall write $SL_2$ for $SL(X^\ast)$ and $B_2$ for its Borel subgroup.

\begin{Def} For any smooth admissible representation $\pi$
  its   Fourier Jacobi module of $\pi$ over $SL_2\times O(V_3)$,
  is defined by
  $${FJ_\psi(\pi)=\Hom_{H}(\omega_{\psi,q}, \pi_{Z,\psi})}.$$ 
\end{Def}

  We record several  properties of the Fourier-Jacobi functor $FJ_\psi$.
\begin{Prop}\label{FJ:prop}
  \begin{enumerate}
    \item
The functor $FJ_\psi$ is  exact.
\item
  $FJ_\psi(\pi)=0$ for an irreducible $\pi$  if and only if
  $\pi$ is one-dimensional. 
 \item 
$FJ_\psi(I_{2n}(\chi,s))=\Ind^{SL_2}_{B_2}(\chi\chi_K,s)\otimes 1.$
\end{enumerate}
\end{Prop}
\begin{proof}
  It is shown in \cite{Weissman}, page $282$
  that $\omega_{\psi,q}$ is a projective object in the category of smooth
  representations of $H$. 
  The Fourier-Jacobi functor is a composition of two exact functors and hence
  exact.

  The representation  $FJ_\psi(\pi)$ is zero
  if and only if $Z$ acts trivially on $\pi$.
  Hence the normal group generated by all the conjugates of $Z$ also acts
  trivially. Since the quotient by any normal non-central subgroup of $O(V_1)$
is finite abelian, the claim follows. 

The proof of part 3 is computational, completely analogous
to the proof of Theorem $4.3.1$ in \cite{Weissman}.
\end{proof}  

\begin{proof}[Proof of Proposition  \ref{red:ind}]
  By Proposition \ref{FJ:prop}, the representation $I_{2n}(\chi,s)$
  is reducible if and only if it either contains
  a one-dimensional constituent, i.e. $\chi=\chi_0, s=\pm(n-1),$
  or its Fourier-Jacobi functor $\Ind^{SL_2}_{B_2}(\chi\chi_K,s)\otimes 1$
  is reducible, i.e. $\chi=\chi_K,s=\pm 1$.
  \end{proof}

Let us investigate the points of reducibility in more detail.
\begin{itemize}
\item
 
  $FJ_\psi(I_{2n}(\chi_0,\pm(n-1)))$ is irreducible and hence
  $I_{2n}(\chi_0,\pm(n-1))$ has length $2$.
  The trivial representation is a unique
  submodule of $I_{2n}(\chi_0,1-n)$ and is a  unique quotient
  of $I_{2n}(\chi_0,n-1)$. 

\item  
 The length of $I_{2n}(\chi_K,\pm 1)$
  equals the length of $FJ_\psi(I_{2n}(\chi_K,\pm 1))$ which is two.
  Since $FJ_\psi(I_{2n}(\chi_K,-1))$ has the trivial representation
  as its unique submodule,  $I_{2n}(\chi_K,-1)$
  has unique submodule, which we denote $\Pi_1$
  and $FJ_\psi(\Pi_1)$  is the trivial representation. This implies
by  \cite{GanSavin}, Proposition $3.7$ that the character expansion
of $\Pi_1$ is $c_1 \hat \mu_{\mO_{\min}}+c_0$, i.e. $\Pi_1$ is minimal.

The representation $\Pi_1^\vee$ is a unique quotient of $I_{2n}(\chi_0,1)$.
On the other hand the intertwining operator $\mM(\chi_K,1)$
defines a non-zero map in $\Hom_{G_1}(I_{2n}(\chi_0,1), I_{2n}(\chi_0,-1))$
whose image is $\Pi_1$. Hence $\Pi_1$ is self-contragredient.  

In the next section  the minimal representation
is embedded in a space of functions on isotropic cone.
The properties of the minimal representation below are used in the proof. 

\begin{Prop}\label{Pi1:prop}
  \begin{enumerate}
\item The representation $\Pi_1$ is unitarizable.
  \item  $(\Pi_1)_{N_1, \Psi}=0$
    for all the anisotropic characters $\Psi$ of $N_1$. 
\end{enumerate}
\end{Prop}

\begin{proof}
  \begin{enumerate}
    \iffalse
  Let $\chi=\chi_K|\cdot|^{-1}$.
The representation $\bar \Pi_1$ is a submodule of
$I_{2n}(\bar \chi)= I_{2n}(\chi)$ and hence $\bar \Pi_1\simeq \Pi_1$.
Besides $\Pi_1^\vee$ is a quotient of $I_{2n}(\chi)^\vee\simeq I_{2n}(\chi^{-1})$
and so is $\Pi_1$ as the image of $\mM(\chi^{-1})$.
Since $I_{2n}(\chi)$ has length two,
it follows that $\Pi_1$ is self contragredient.
So $\bar \Pi_1\simeq \Pi_1^\vee$, i.e. $\Pi_1$ admits a non-degenerate
Hermitian form. Let us show that it is positive definite.
\fi
\item
For $s\ge 0$ define a non-zero $G_1$-invariant Hermitian form $H(s)$
on $I_{2n}(\chi_K,s)$
by $H(s)(f)=
\int\limits_{Q\backslash G} f(x)\overline{\mM(\chi_K,s)(f)(x)}dx$

For $s=0$ the representation $I_{2n}(\chi_K,0)$ is unitary irreducible.
Hence it admits unique up to constant $G_1$-invariant Hermitian form
which is positive definite. Hence $H(0)$ is definite. By multiplying
by $(-1)$, if needed, we can assume that $H(0)$ is positive definite.

For any  $0\le s <1$ the representation $I_{2n}(\chi_K,s)$
is irreducible, so $H(s)$ is positive definite and non-degenerate.
For $s=1$ the form is positive semi-definite
and is degenerate. Indeed, the image of $\mM(\chi_K,1)$ is $\Pi_1$
and hence the kernel of $\mM(\chi_K,1)$ is the radical of $H(1)$.
Thus $H(1)$ reduces to a non-degenerate positive definite form
on $\Pi_1$ and so  $\Pi_1$ is unitarizable. 
\item
In \cite{MoeglinWaldspurger} to any coadjoint nilpotent orbit $\mO$
is associated a degenerate Whittaker model $W_{\mO}(\Pi_1)$
and the set of maximal elements in ${\{\mO: W_\mO(\Pi_1)\neq 0\}}$
coincides with the wave front of $\Pi_1$.
%For $\mO=\mO_{min}$ one has $W_{\mO}(\Pi_1)=FJ_\psi(\Pi_1)$
%is one-dimensional. Hence the coefficient $c_1=1$. 
If $(\Pi_1)_{N_1,\Psi}\neq 0$ for an anisotropic character $\Psi$ then
$W_{\mO}(\Pi_1)\neq 0$ for the Richardson nilpotent orbit
$\mO$ corresponding to $N_1$. By \cite{MoeglinWaldspurger}
this contradicts the minimality of $\Pi_1$.
\end{enumerate}
  \end{proof}
\end{itemize}

\subsection{The space of functions on the cone}
From now on for any character $\chi$ we denote by $I_{2n}(\chi)$
the normalized induced representation $\Ind^{G_1}_{Q_1}\chi$.
This representation and  in particular, the minimal representation $\Pi_1$ can be realized
on a space of functions on the cone.

To simplify notation we shall write $\mS^\infty$ (resp. $\mS_c$)
for the space of  smooth functions (resp.
 smooth functions of compact support) on $C_0$. 

The group $G_1$ acts transitively on $C_0$ and the derived group
$[Q_1,Q_1]=O(V_2)\cdot N_1$ is the stabilizer subgroup of $e_1^\ast$.

The space $\mS_c$ with the action of $GL_1\times G_1$ defined in
\ref{intro:Q:action} is isomorphic as representation of $GL_1\times G_1$  
to $\ind^{GL_1\times G_1}_{GL_1\times Q_1} \mS_c(F^\times)$,
where the action of $GL_1\times M_1$ is given by
$$\left\{
\begin{array}{ll}
  a\cdot \varphi(x)=\chi_K(a)|a|^{n-1}\varphi(ax)& a\in GL_1\\
  (b,h)\cdot \varphi(x)=\varphi(b^{-1}x) & (b,h)\in M_1
\end{array}\right.
  $$
and extended trivially to $GL_1\times Q_1$. 

%This implies
%$(\mS_c)_{GL_1,\chi^{-1}}\simeq \chi^{-1}\otimes I_{2n}(\chi_K\chi)$
%for any character $\chi$ of $GL_1$. 

\begin{Def}\label{Mellin:p} For any character $\chi$ of $GL_1$
the  Mellin transform,
${p_{\chi}:\mS_c\rightarrow \mS^\infty}$ is  defined by 
$$p_{\chi}(f)(w)=\int\limits_{GL_1} (a\cdot f)(w)\chi(a)\, da.$$
\end{Def} 

The image of $p_\chi$, denoted by $\mS(\chi),$ is isomorphic
 to $\chi^{-1}\otimes I_{2n}(\chi\chi_K)$ as $GL_1\times G_1$ module. 

There exists  unique, up to a constant,
$G_1$-invariant measure $\omega$ on $C_0.$
Our choice of $\omega$ is fixed in \ref{def:measure}. 

\begin{Prop}\label{mult:one:Phi1}
  \begin{enumerate}
  \item  For any $\chi$ the space
    \begin{equation}\label{hom:onedim}
      \Hom_{GL_1\times G_1}(\mS_c,\chi^{-1}\otimes I_{2n}((\chi\chi_K)^{-1})
      \end{equation}
      is one-dimensional  and is spanned by
$\mM(\chi\chi_K)\circ p_\chi$. 
    \item For $\chi=\chi_K|\cdot|^{n-1}$ the image of any non-zero      element in \ref{hom:onedim} is the one-dimensional space of constant functions,
      isomorphic to       ${\chi_K|\cdot|^{1-n}\otimes 1}$  as $M$-module.
      The map is given by  
\begin{equation}\label{int:cone}
  f\mapsto I_C(f)=\int\limits_{C} f(v)\omega(v)
  \end{equation}
\item For $\chi=|\cdot|$,   the image of any non-zero
  element in \ref{hom:onedim} is
 ${|\cdot|^{-1}\otimes \Pi_1}$ as  $M$-module.
\end{enumerate}
  \end{Prop}

\iffalse
\begin{Prop}\label{mult:one:Phi1}
  \begin{enumerate}
    \item 
      The space
      $$\Hom_{GL_1\times G_1}(\mS_c,
      \chi_K|\cdot|^{1-n}\otimes \mS(\chi_K|\cdot|^{1-n}))$$
  is one-dimensional  and is generated by the functional 
\begin{equation}\label{int:cone}
  f\mapsto I_C(f)=\int\limits_{C} f(v)\omega(v)
  \end{equation}
  The image of any non-zero homomorphism is one-dimensional representation
  $\chi_K|\cdot|^{1-n}\otimes 1$  as $M$-module. 
\item 
 The space $$\Hom_{GL_1\times G_1}(\mS_c, |\cdot|^{-1}\otimes \mS(|\cdot|^{-1}))$$
 is one dimensional and the image of any non-zero homomorphism is
 ${|\cdot|^{-1}\otimes \Pi_1}$ as  $M$-module.
\end{enumerate}
  \end{Prop}
\fi
\begin{proof}
For an arbitrary character $\chi$ holds:
  % $$\Hom_{GL_1\times G_1}(\mS_c, \chi^{-1}\otimes \mS(\chi))=
$$\Hom_{GL_1\times G_1}(\mS_c, \chi^{-1}\otimes I_{2n}((\chi_K\chi)^{-1})=$$
$$\Hom_{GL_1\times G_1}((\mS_c)_{GL_1,\chi^{-1}},
\chi^{-1}\otimes I_{2n}((\chi_K\chi)^{-1}))=
\Hom_{G_1}(I_{2n}(\chi_K \chi), I_{2n}((\chi_K\chi)^{-1}))$$
that is one-dimensional, spanned by $\mM(\chi\chi_K).$
The image of $\mM(\chi\chi_K)$ is unique irreducible submodule of
$I_{2n}((\chi_K\chi)^{-1})$. 
  \end{proof}

\subsection{The Schrodinger model}\label {subsec:schrodinger}.

The results of the previous subsection, applied to the group $G$,
show that  the degenerate
principal series ${I_{2n+2}(\chi_K|\cdot|^{-1})}$ contains a unique
irreducible submodule $\Pi$, which is a minimal representation of $G$.

In this section we describe an embedding of $\Pi$ in $\mS^\infty$ and 
prove Theorem \ref{savin:intro}.
The image $\mS$ of the embedding is often referred to as
Schrodinger model of $\Pi$. The proof in the split case
appears in \cite{Savin} and a  more detailed proof of the same theorem in
\cite{SavinWoodbury}. The proof is essentially  the  same
for quasi-split groups and we write it below, filling in some details
for the reader's convenience. 
   
  \begin{proof}[Proof of Theorem \ref{savin:intro}]
    \begin{enumerate}
      \item
    Let $\bar I_{2n+2}(\chi)$ denote the degenerate principal series
    induced from  the character $\chi$ of the opposite parabolic
    subgroup $\bar Q$. It contains a $Q$-subspace $\bar I^{0}_{2n+2}(\chi)$
    consisting of functions, whose support is contained
    in the open Bruhat cell $\bar Q N$.

    Since $\bar QN$ is open dense in $G$, the restriction to $N$
    defines a $Q$ equivariant
    embedding  $\bar I_{2n+2}(\chi_K|\cdot|)\hookrightarrow \mS^\infty(N)$
    that maps $\bar I^{0}_{2n+2}(\chi_K|\cdot|)$ onto $\mS_c(N).$
    
    The action of $Q$ on $\mS^\infty(N)$ is given by   
    \begin{equation}
      \left\{\begin{array}{ll}
      n\cdot  f(x)=f(xn)& x,n\in N\\
      m\cdot f(x)=\chi_K(m)|m|^{1-n}f(m^{-1}xm) & m\in M
\end{array}\right..
      \end{equation}

    Define a   non-degenerate $G$-invariant Hermitian pairing
    $$\la, \cdot, \cdot\ra:\bar I_{2n+2}(\chi_K|\cdot|)\times
    \bar I_{2n+2}(\chi_K|\cdot|^{-1})\rightarrow \C $$
    $$\la f,g\ra=\int\limits_{\bar Q\backslash G} f(x)\overline{g(x)}dx=
    \int\limits_{N} f(n)\overline{g(n)}dn.$$

    The representation $\Pi$ is embedded in
    $\bar I_{2n+2}(\chi_K|\cdot|)$ as a unique submodule. 
    For any $0\neq f\in \Pi\subset \mS^\infty(N)$ there is
    $g\in \mS_c(N)$ such that $\la f,g\ra\neq 0,$
    i.e. the pairing $\la \cdot,\cdot\ra$ restricted to
    $\Pi\times \bar I_{2n+2}(\chi_K|\cdot|^{-1})$ is left non-degenerate.
  
    Denote by $\mF_{\psi,N}:\mS_c(N)\rightarrow\mS_c(\bar N)$
    the Fourier transform with respect to the Killing form. 
    Recall that $C\subset \bar N$. Define  the $Q$-submodules
    $\mS''\subset \mS'\subset \mS_c(N):$
    $$ \mS''= \{f: \mF_{\psi,N}(f)|_C=0\}\subset
    \mS'=\{f: \mF_{\psi,N}(f)(0)=0\}\subset \mS_c(N)$$
       
    Then $\mF_{\psi,N}: \mS'/\mS''\simeq \mS_c,$
    where the action
    of $Q$ on $\mS_c$ by the transport of structure is given by:

    \begin{equation}
      \left\{\begin{array}{ll}
      n\cdot  f(x)=\psi(-\la n,x\ra)f(x)&   n\in N\\
      m\cdot f(x)=\chi_K(m)|m|^{n-1}f(m^{-1}x m)& m\in M
\end{array}\right.
      \end{equation}
Here $x\in C\subset \bar N$.

    Our goal is to show that the pairing $\la \cdot, \cdot\ra$ reduces
    to a left non-degenerate  pairing $\Pi\times \mS_c$.
    For this we need to show
    \begin{Lem}
    \begin{enumerate}
    \item The pairing restricted to $\Pi\times \mS'$ is left non-degenerate. 
    \item The pairing restricted to $\Pi\times \mS''$ is zero.
    \end{enumerate}
    \end{Lem}
\begin{proof}    
  \begin{enumerate}
    \item
      Assume that
    $0\neq f\in \Pi$, such that $\la f,g\ra =0$ for all $g\in
    \mS'.$
    The space $\mS'$ consists of functions $g\in \mS_c(N)$ such that
    $\int\limits_N g(n)dn=0$. 
    For any $h\in \mS_c(N)$ and $n\in N$,
    the function $nh-h$ belongs to  $\mS'$. Hence for any $n\in N$
    holds 
    $$0=\la f,nh-h\ra=\la n^{-1}f-f,h\ra, \quad \forall h\in \mS_c(N).$$
    In particular, the vector $f\in \Pi$ is $N$-invariant.
    However, $\Pi$ is unitarizable by \ref{Pi1:prop}, part $1$,
    and hence by the result of Howe and Moore \cite{HoweMoore} does not contain
    $N$-fixed vectors. This is a contradiction.

\item
    Assume that there exists $0\neq g\in \mS''$, i.e. $\mF_{N,\psi}(g)|_{C}=0$
    and $f\in \Pi$ such that $\la f, g\ra\neq 0$.

    Choose an open compact subgroup $\mN_1\subset N$
    such that the supports of both $g$ and $\mF_{\psi,N}(g)$ 
    are contained in $\mN_1$. One has  
$$\la f,g\ra =\int\limits_N f(n)\overline{g(n)}\, dn =
\int\limits_{\mN_1} f(n)\overline{g(n)}\,dn=$$
$$\int\limits_{\mN_1} f(n)
\int\limits_{\mN_1} \overline{\mF_{\psi,N}(g)(u)}\Psi_{u}(n)\, du\, dn=
\int\limits_{\mN_1}  l_{\mN_1,u}(f)\overline{\mF_{\psi,N}(g)(u)}du,$$
where 
$$l_{\mN_1,u}(f)=\int\limits_{\mN_1} f(n)\Psi_{u}(n) dn.$$

The function $\mF_{\psi,N}(g)$ is supported on anisotropic vectors.
Hence there exists an anisotropic vector $u$ such that $l_{\mN_1,u}(f)\neq 0$.

On the other hand  $\Pi_{N,\Psi_{-u}}=0$ by \ref{Pi1:prop}, part $2$. Hence 
we can write  $f=\sum_{i=1}^k n_i f_i-\Psi_{-u}(n_i) f_i$
for $n_i\in N, f_i\in \Pi$. Without loss of generality we can assume
that $n_i\in \mN_1$ for all $i$  and so  $l_{\mN_1,u}(f)=0$. This is a contradiction.
\end{enumerate}
  \end{proof}

Hence the pairing reduces to the left non-degenerate
$Q$-invariant pairing on $\Pi\times \mS_c$. This implies that
$\Pi$ is embedded into a $Q$-smooth dual of $\mS_c$
that can be described as a space $\mS$ of smooth functions,
not necessarily of compact support, with the  action of $Q$  given by 
    \begin{equation}\label{Q:action:on:S}
      \left\{\begin{array}{ll}
      n\cdot  f(x)=\psi(-\la n,x\ra)f(x)& n\in N\\
      m\cdot f(x)=\chi_K(m)|m|^{n-1}f(m^{-1}xm) & m\in M
\end{array}\right.
    \end{equation}    
       This completes part one.
\item    
Let us show that $\mS$ contains $\mS_c$.  Since  $\Pi$ is a quotient of
   $\bar I_{2n+2}(\chi_K|\cdot|^{-1}),$ the $G$-invariant pairing
   $\la\cdot,\cdot,\ra$  gives rise to the
   non-degenerate pairing on $\Pi\times \Pi$. In particular $\mS''$
   is contained in the irreducible submodule of  $\bar I_{2n+2}(\chi_K|\cdot|^{-1}).$
   
The composition of  $Q$-equivariant maps
$$\mS_c\simeq \mS'/\mS''\hookrightarrow \bar
I_{2n+2}(\chi_K|\cdot|^{-1})/\mS''\rightarrow \Pi$$
is not zero. The space $\mS_c$ is an irreducible $Q$-module,
hence $\mS$ contains $\mS_c$. 

Finally, any vector $f\in \mS$ is fixed by a compact subgroup of $N.$
The formula of action of $N$ implies that $f$ is  of bounded support.       
\item
   The space  $\mS[N]=\{n\cdot f -f|f\in \mS\}$ is contained in
    $\mS_c$. On the other hand $(\mS_c)_N=0,$ i.e. $\mS_c=\mS_c[N]$. 
    Hence $\mS_c[N]\subseteq \mS[N]\subseteq \mS_c=\mS_c[N].$
Hence $\mS[N]=\mS_c$ and  so $j_N(\mS)\simeq \mS/\mS_c$,  
   as required.
\end{enumerate}
    \end{proof}  

  \iffalse
  The space $L^2(C_0)$ is an irreducible representation of $Q$. Hence
  the unique (up to constant) irreducible $Q$-invariant inner product of
  $\mS_c(C_0)$ is given by  the standard inner product on $C_0$
  with respect to the $G_1$-invariant measure $\omega$.
  Since $\Pi$ is unitarizable and $\mS_c\subset \mS$ it follows that
  $\mS$  is contained  $L^2(C_0)$ and
  is obviously dense there. Hence $\mS$ is the space of smooth functions
  of the unitary minimal representation $\hat \Pi,$ realized on $L^2(C)$.

Since $Q$ is a maximal parabolic group, the group $G$
is generated by $Q$ and any element not contained in $Q$.
An explicit  formula for the action of such element will
complete the action of the group.
\fi

Let $r\in O(V)$ be the involution defined in \ref{def:r}.
The following Lemma reduces the description of
the space of germs $\{[f]_0, f\in \mS\}$ to the
description of the space $\{[\Pi(r)f]_0, f\in \mS_c\}$. 

\begin{Lem}\label{S=S0+Phi(S0)}
  $\mS=\mS_c+\Pi(r)(\mS_c).$
\end{Lem}

\begin{proof} The RHS is obviously included in LHS, since the space $\mS$ is
  $G$-invariant.
  Since $\mS$ is irreducible, it is enough to show that RHS is
  $G$-invariant too.
  It is invariant under the action of $G_1$ which preserves $\mS_c$
  and commutes with $r$.
  It is  invariant under $r$, since $r$ is involution.
  It is also invariant under the action of $N$, since $N$ preserves
  $\mS_c$ and acts trivially  on $j_N(\mS)=\mS/\mS_c$.
  The groups $G_1, N$ together with the element $r$
  generate the group $G$ and hence RHS is $G$-invariant.
  \end{proof}

\section{ Radon transform $\mR$ and Fourier transform $\Phi$  }\label{sec:radon}

In this section we define the  Radon transform on a cone
and study its asymptotics. This enables us to define the Fourier
transform on a cone. 

\subsection{The Fourier transform on a quadratic  space and Weil index}
\label{subsec:FT}
Let  $(U,q)$ be   a non-degenerate quadratic space  of
dimension $2m$ with the associated bilinear form $\la \cdot, \cdot\ra$.
There exists  unique Haar measure $du$ on $U$,
called self-dual with respect to $\psi,q$ such that
the Fourier transform $\mF_{\psi,U}$ on $\mS_c(U)$ defined by
$$\mF_{\psi,U}(f)(v)=\int\limits_{U}f(u) \psi(\la u,v\ra) du $$
is unitary and satisfies
$\mF_{\psi,U}\circ \mF_{\psi,U}(f)(v)=f(-v)$ for all $f\in \mS_c(U)$.
We shall always use this measure on quadratic spaces without further
notice. 

For any $t\in F^\times$ let $\psi_t$ be an additive character
defined by $\psi_t(x)=\psi(tx)$. It is easy to see that
\begin{equation}\label{Fourier:psi:t}
  \mF_{\psi_t,U}(f)(v)=|t|^m\mF_{\psi,U}(f)(tv).
  \end{equation}

\vskip 5pt
If $U=F$ is one-dimensional and $q(x)=x^2$ we omit $U$
and write $\mF_\psi$ for the Fourier transform of one variable.
The self-dual measure $dx$ on $F$ gives rise to the Haar measure
$d^\times x=\frac{dx}{|x|}$ on $F^\times$.

\vskip 5pt
Let $K$ be a quadratic algebra over $F$,
with the Galois involution $x\mapsto \bar x,$
a norm $\N$ and the trace $\T$. The bilinear form $\la x,y\ra=\T(x\bar y)$
gives rise to  the quadratic form $q(x)=\N(x)$.

According to Weil, \cite{Weil}
there exists a constant $\gamma(\chi_K, \psi),$ which
is  a fourth root of unity satisfying
\begin{equation}\label{weil:K}
\int\limits_K \mF_{\psi,K}(f)(x) \psi(\N(x))dx=
\gamma(\chi_K,\psi)\int\limits_K f(x)\psi(-\N(x)) dx.
\end{equation}
Moreover, it holds for all $t\in F^\times$
\begin{itemize}
  \item
$\gamma(\chi_K,\psi_t)=\chi_K(t)\gamma(\chi_K,\psi),$
\item
$\gamma(\chi_K, \psi)=1$ if $K$ is split. 
\end{itemize}

\iffalse
For $Re(s)>0$ define a function $\varphi_s$ such that
$\varphi_s(x)=|x|_K^s$ for $|x|\le 1$ and zero otherwise.
It is continuous function, but not locally constant at zero.
It is easy to see that the identity \ref{weil:K} still holds for
$f=\varphi_s$. We denote by $\mS^s_c(K)$ the space 
$\mS_c(K)+\Span_\C \{\varphi_s\}$.
\fi

Assume that the space $U=\bH^{m-1}\oplus K, q=q^{m-1}_{\bH}\oplus \N.$
For any $f\in \mS_c(U)$ holds
\begin{equation}\label{weil:U}
\int\limits_U \mF_{\psi,U}(f)(u) \psi(q(u))du=
\gamma(\chi_K,\psi)\int\limits_U f(u)\psi(-q(u)) du.
\end{equation}

\begin{Lem} For any $t\in F^\times$ holds 
\begin{equation}\label{weil:psi:t}
\int\limits_U \mF_{\psi,U}(f)(u) \psi(tq(u))du=
|t|^{-m}\gamma(\chi_K,\psi)\chi_K(t)\int\limits_U f(u)\psi(-t^{-1}q(u)) du.
\end{equation}
\end{Lem}
\begin{proof}
Plug  the character $\psi_{t^{-1}}$ in \ref{weil:U} to get
$$\int\limits_U \mF_{\psi_{t^{-1}},U}(f)(u) \psi(t^{-1}q(u))du=
\gamma(\chi_K,\psi_{t^{-1}})\int\limits_U f(u)\psi(-t^{-1}q(u)) du.$$
Note that the measure can be still taken the self-dual with respect to
$\psi,q$ on both sides.
Using \ref{Fourier:psi:t}, the property of the Weil index mentioned above 
and the change of variables $t^{-1}u\mapsto u$ we obtain the required.
%$$|t|^{m}\int\limits_U \mF_{\psi,q}(f)(t^{-1}u) \psi(tq(t^{-1}u))dt^{-1}u=
%\gamma(\chi_K,\psi)\chi_K(t)\int\limits f(u)\psi(-t^{-1}q(u)) du$$

\end{proof}
%REPLACE
\subsection{Gelfand-Leray forms and measures}\label{gelfand:leray}
Let $Y$ be a smooth algebraic variety over $F$
and $\varphi:Y\rightarrow F$ be an algebraic map. 
Let $Y^{sm}\subset Y$ denote the set of smooth points
of $\varphi$. We write $Y(t)$ for the fiber $\varphi^{-1}(t)$
and $Y(t)^{sm}=Y(t)\cap Y^{sm}$ for the smooth fiber.

Let $\eta$ be a smooth differential top-form on $Y$.
Locally on $Y^{sm}$ there exist smooth forms $\tilde \eta$
such that $\tilde \eta\wedge \varphi^\ast(dx)=\eta$.
Since the restriction of a form $\tilde \eta$ to fibers
does not depend on the choice of
$\tilde\eta$ we obtain a family $\eta_t$ of {\it Gelfand-Leray} 
top forms  on $Y(t)^{sm}$. We  denote by $|\eta_t|$ the associated
measures on  $Y(t)^{sm}$. It is easy to see that 
$$\int\limits_F \int\limits_{Y(t)^{sm}} f(y) |\eta_t(y)| dt=
\int\limits_{Y^{sm}} f(y) |\eta(y)|$$
for any $f\in\mS_c(Y^{sm})$.

We shall use this construction in the following cases:
\begin{enumerate}
  \item
  For a non-degenerate quadratic space $(U,q)$ let $Y=U$,
  $\varphi=q$ and $U^{sm}=U-\{0\}$. Then $U(t)^{sm}=U(t)$
  unless $t=0$ and $U(0)^{sm}=U(0)-\{0\}$. 
  Choose $\eta=du$, the self-dual measure with respect to $(\psi,q)$.
  The form $\eta_0$ on $U(0)^{sm}$ gives rise
  to a $O(U)$-invariant measure $|\eta_0|$. 

  Define a family of distributions $\delta_{U}(t)$ on $\mS_c(U(t)^{sm})$  by
\begin{equation}\label{def:delta}
  \delta_{U}(t)(f)=\int\limits_{U(t)^{sm}} f(u) |\eta_t(u)|,
  \end{equation}
  so that by Fubini theorem for any $f\in \mS_c(U^{sm})$ holds
  $$\int\limits_F \delta_{U}(t)(f) dt=\int\limits_{U} f(u) du$$ 

  We show in Lemma \label{delta:convergence}  that
  the integral \ref{def:delta} converges for all $f\in \mS_c(U(t)),$
or equivalently for all $f\in \mS_c(U(0))$.

  \begin{Def}\label{def:measure}
    Take $(U,q)=(V_1,q)$, so that $U(0)^{sm}=C_0$.
    The construction above gives rise 
    to the  $O(V_1)$ invariant measure on $C_0$ that we denote by $\omega$.
\end{Def}

\item
  Let $Y=C_0\subset V_1$. For each $w\in C_0$  define a map
  $$\varphi_w:C_0\rightarrow F, \quad \varphi_w(v)=\la v,w\ra$$
  and let $Y_w(t)$ be the fiber over $t$. 
  The set of smooth points $C_w^{sm}$ of $\varphi_w$ is
  $C_0-L_w$, where $L_w=Span\{w\}$. Starting with
  the measure $\omega$ on $C_0$,
  we obtain for each $t$ the measures  $\omega_{w,t}$ on
  the smooth fibers $Y^{sm}_w(t).$  

For any $w\in C_0$ consider a family of distributions
$\mR_w(t)$ on $\mS_c(C_w^{sm})$, defined by 
  \begin{equation}\label{def:radon}
    \mR_w(t)(f)=\int\limits_{Y^{sm}_w(t)} f(v) \omega_{w,t}(v), \quad
    f\in\mS_c(C_w^{sm}).
    \end{equation}

  We  show in Lemma
  \ref{lemma:radon:convergence} that the integral \ref{def:radon}
  converges for all $f\in\mS_c$. Its convergence depends
  on the convergence of the integral \ref{def:delta} on a subspace
  $V_w$ of $V_1$ of codimension $2$.
  
Once this is proven,  we obtain a family of operators
$$\mR(t):\mS_c\rightarrow \mS^\infty,
\quad \mR(t)(f)(w)=\mR_w(t)(f)\quad  \forall w\in C_0.  $$ 
These operators are called Radon transforms on the cone. 
\end{enumerate}

\subsection{The distributions $\delta_{V_2}(t)$.}
\begin{Lem}\label{delta:convergence} Let $\dim V_2=2n-2\ge 4$.
  The integral \ref{def:delta} defining $\delta_{V_2}(t)(f)$
  converges for any $f\in \mS_c(V_2)$.
\end{Lem}
\begin{proof}
  For $t\neq 0$ the restriction of $f$ to $V_2(t)=V_2(t)^{sm}$
  is of compact support and hence the integral converges.
  For $t=0$ it is enough to check that the integral \ref{def:delta} converges
  for the function $f=1_{V_2(\mO_F)}$.
  We write $1_{V_2(\mO_F)}=\sum_{k=0}^\infty r_k$, where
  $r_k \in\mS_c(V_2^{sm})$ is the characteristic function of
  ${\{v\in V_2(\mO_F): |v|=q^{-k}\}}$.
By homogeneous property of $\delta_V(0)$ we have  
$$\delta_{V_2}(0)(1_{V_2(\mO_F)})=
\sum_{k=0}^\infty q^{-(2n-4)k}\delta_{V_2}(0)(r_0)$$
  which converges for $n>2$. 

\end{proof}

The following proposition plays an important role in the paper.
The second part, which follows from the first,
describes the asymptotic behavior of $\delta_{V_2}(s)$ for small $s$.
It exhibits $\delta_{V_2}(s)(f),$ which is an integral of
$f$ over the hyperboloid $V_2(s),$ as a sum of the
integral of $f$  over the cone $V_2(0)$
and the error term, which depends on the value  $f(0).$
This part is used  to describe
the asymptotic behavior of the Radon transform $\mR(t)$ for small $t$
in section \ref{subsec:Radon}, which in turn allows to define the
Fourier transform on a cone. It is also used in the paper
\cite{GurevichKazhdanBAS}.

The first part will be used in Theorem \ref{Pi(r)=Phi}.
In fact, the decomposition of $j_N(\mS)$
as a sum of two irreducible representations
is a result of the decomposition of $\delta_{V_2}(s)$ in part $1$.
  
\begin{Prop} Let $\dim V_2=2n-2\ge 4$ and  $f\in \mS_c(V_2)$.
  \begin{enumerate}
    \item
      $$\delta_{V_2}(s)(\mF_{\psi,V_2}(f))=\delta_{V_2}(0)(\mF_{\psi,V_2}(f))+
      \int\limits_{V_2} f(v) H_s(v)dv,$$
where
\begin{equation}\label{H:def}
      H_s(v)=
  \gamma(\chi_K,\psi)\chi_K(s)|s|^{n-2}
  \int\limits_F  \psi(q(v)t s) \chi_K(-t)(\psi(t^{-1})-1)  |t|^{n-2} d^\times t
  \end{equation}
\item For any $f$ there exists $\epsilon>0$ such that for $|s|<\epsilon$
  holds
\begin{equation}\label{delta:difference}
  \delta_{V_2}(s)(f)=\delta_{V_2}(0)(f)+ c_{\psi,q}\chi_K(s)|s|^{n-2}f(0),
  \end{equation}
 where
   $$c_{q,\psi}=\gamma(\chi_K,\psi)
   \int\limits_F (\psi(t^{-1})-1) \chi_K(-t)|t|^{n-2} d^\times t.$$
\end{enumerate}
\end{Prop}

\begin{proof}
  It is easy to see that part $1$ implies part $2$. Indeed,
  for fixed $f\in S_c(V_2)$ there exists $\epsilon $ such that $|s|<\epsilon$
  implies $\psi(sq(v)t)=1$ for all $v\in supp(f)$
  and $t\in supp(\psi(t^{-1})-1).$  Hence, by part $(1)$,
  \begin{multline}
    \delta_{V_2}(s)(\mF_{\psi,V_2}(f))=\delta_{V_2}(0)(\mF_{\psi,V_2}(f))+
    c_{\psi,q}\chi_K(s)|s|^{n-2}\cdot
    \int\limits_{V} f(v)dv=\\
    \delta_{V_2}(0)(\mF_{\psi,V_2}(f))+ c_{\psi,q}\chi_K(s)|s|^{n-2}\mF_{\psi,V_2}(f)(0).
    \end{multline}

  Since the operator $\mF_{\psi,V_2}\in Aut(\mS_c(V))$ is surjective
  the equation (\ref{delta:difference}) follows. 

Let us prove part $(1)$.
For $t\in F$ the function
$\psi_t\circ q$ defines a distribution $q(\cdot, t)$ on $\mS_c(V_2)$ by
$$q(f,t)=\int\limits_{V_2} f(v) \psi(tq(v)) dv.$$

For any $f\in \mS_c(V_2)$ the  function $t\mapsto q(f,t)$
is locally constant and belongs to $L^1(F)$. Indeed,  
 it follows from (\ref{weil:psi:t}) that 
 $$q(\mF_{\psi,V_2}(f),t)=
 |t|^{1-n}\chi_K(t) \gamma(\chi_K,\psi) \cdot q(f,-t^{-1}).$$
 Thus for large $|t|$ holds $|q(\mF_{\psi,V_2}(f),t)|=C|t|^{1-n}$ 
 and hence $q(\mF_{\psi,V_2}(f),\cdot)$ belongs to $L^1(F)$.
 By surjectivity of $\mF_{\psi,V_2}\in Aut(\mS_c(V_2))$ one has
 $q(f,\cdot)\in L^1(F)$ for any $f$.

  By Fubini theorem for any $f\in \mS_c(V_2)$ holds
  $$q(f,t)=\int\limits_{V_2} f(v)\psi(tq(v))dv=
  \int\limits_{F}\int\limits_{V_2(s)} f(v)|\eta_s(v)| \psi(ts) ds=$$
  $$\int\limits_F \delta_{V_2}(s)(f) \psi(ts) ds=
  \mF_\psi(\delta_{V_2}(\cdot)(f))(t).$$
  Applying the inverse Fourier transform to both sides we get
 $$\delta_{V_2}(s)(f)=\int\limits_F q(f,t) \psi(-st) dt.$$
Substituting  $\mF_{\psi,V_2}(f)$ instead of $f$ and making
the change of variables $t\mapsto (-st)^{-1}$  this yields
$$\delta_{V_2}(s)(\mF_{\psi,V_2}(f))-\delta_{V_2}(0)(\mF_{\psi,V_2}(f))=
\int\limits_F q(\mF_{\psi,V_2}(f),t) (\psi(-st)-1)dt =$$
$$\gamma(\chi_K,\psi)
  \int\limits_F |t|^{1-n} \chi_K(t)
  \int\limits_{V_2}  f(v)\psi(-\frac{q(v)}{t})(\psi(-st)-1) dv dt=$$
$$\gamma(\chi_K,\psi)\chi_K(s)|s|^{n-2}
\int\limits_F |t|^{2-n}\chi_K(-t)
\int\limits_{V_2} f(v)\psi(\frac{q(v)s}{t})dv\cdot  (\psi(t)-1) d^\times t= $$
$$\gamma(\chi_K,\psi)\chi_K(s)|s|^{n-2}
\int\limits_F 
\int\limits_{V_2} f(v)\psi(q(v)st) dv  (\psi(t^{-1})-1)\chi_K(-t) |t|^{n-2}
  d^\times t. $$
Since  the term $\psi(t^{-1})-1$ vanishes for $|t|$ large
this can be rewritten as 
$$\chi_K(s)|s|^{n-2}\int\limits_{V_2} f(v)
\left(\gamma(\chi_K,\psi)
\int\limits_F
\psi(q(v)t s)(\psi(t^{-1})-1) \chi_K(-t) |t|^{n-2} d^\times t
\right)
dv,$$ 
as required.
\end{proof}

\subsection{The Radon transform $\mR(t)$}\label{subsec:Radon}
Recall that  $C$ is a cone of isotropic vectors in
$(V_1,q)$ with $\dim(V_1)=2n$.

Fix  $w\in C_0$  and choose a vector $w^\ast\in C_0$ such that
$\la w, w^\ast\ra=1$. 
Then  $\bH_w=Span\{w,w^\ast\}$ is a hyperbolic plane  and
let $V_w=\bH_w^\perp$ in $V_1$. Any $v$ in $V_1$ can be written uniquely as 
$$v=tw^\ast+u-sw, \quad u\in V_w; t,s\in F.$$

For any $f\in \mS_c$ define the function $f_{t,s}$ on $V_w(st)$ by
$$f_{t,s}(u)=f(tw^\ast+u-sw).$$

We fix a top form $\eta'_w$ on $V_w$ such that
$dv=dt\wedge ds \wedge \eta'_w.$
This gives rise to the family of Gelfand-Leray forms
$\eta'_{w,r}, r\in F$ on $V^{sm}_w(r).$

\begin{Lem}\label{lemma:radon:convergence}
  For $f\in \mS_c,$ the integrals 
  $$\int\limits_{Y^{sm}_w(t)}f(v)\omega_{w,t}(v), \quad  
  \int\limits_F \delta_{V_w}(st)(f_{s,t}) ds$$
both converge and are equal to each other. 
\end{Lem}

\begin{proof}   Let $f\in \mS_c$. 

  Assume  $t\neq 0.$ The map
  $v\mapsto tw^\ast+v-\frac{q(v)}{t}w$ defines bijection
  $V_w\rightarrow Y^{sm}_w(t)$. The function $f_{t,s}$
  belongs to $\mS_c(V_w)$ and $s\mapsto \delta_{V_w}(st)(f_{t,s})$
  is  locally constant and of compact support.  Hence the integral 
  $$\int\limits_{F}\delta_{V_w}(st)(f_{t,s}) ds$$
 converges. By definition of  $\delta_{V_w}(st)$ it equals
 $$\int\limits_F \int\limits_{V_w(st)} f(tw^\ast +v-sw) |\eta'_{w,st}(v)| ds=
 \int\limits_F \int\limits_{V_w(s)}
f(tw^\ast +v-\frac{s}{t}w) |t|^{-1} |\eta'_{w,s}(v)|  ds=$$
$$\int\limits_{V_w} f(tw^\ast +v-\frac{q(v)}{t}w) |t|^{-1} |\eta'_w(v)|=
\int\limits_{Y^{sm}_w(t)} f(v) \omega_{w,t}(v)=\mR_w(t)(f)$$
as required.

Assume $t=0$. The map $(v,s)\mapsto v-sw$ define
a bijection $V_w(0)\times F\rightarrow Y_w(0)$
which restricts to the bijection between smooth points,
i.e. $V_w(0)^{sm}\times F$ is mapped to $Y_w^{sm}(0).$

The function $f_{0,s}$ belongs to $\mS_c(V_w).$
So $\delta_{V_w}(0)(f_{0,s})$ is well-defined and
is a locally constant function of $s$ of compact support.

$$\int\limits_{F}\delta_{V_w}(0)(f_{0,s}) ds=
\int\limits_F \int\limits_{V_w(0)^{sm}} f(v-sw) |\eta'_{w,0}| ds=$$
$$\int\limits_{Y_w^{sm}(0)} f(u) |\omega_{w,0}(u)|=\mR_w(0)(f)$$
as required. 
\end{proof}

Now the Radon operators $\mR(t)$ on $\mS_c$ are well-defined
and the integration on $Y^{sm}_w(t)$ can be replaced
by the integration over $Y_w(t)$. The group $G_1$
acts on the set of fibers, such that
$Y_w(t)g=Y_{wg}(t)$ and $g(\omega_{w,t})=\omega_{wg,t}$.

We shall note two obvious, but useful properties
\begin{itemize}
  \item
By Fubini theorem for any $f\in \mS_c$ and $w\in C_0$ holds
 $$\int\limits_F \mR(t)(f)(w) dt=\int\limits_{C_0} f(v) \omega(v).$$ 
\item
Since $Y_{xw}(xt)=Y_w(t)$ for all $x\in F^\times$
 it holds
 \begin{equation}\label{Radon:homog}
   \mR(xt)(f)(xw)=|x|^{-1}\mR(t)(f)(w).
 \end{equation}
\end{itemize}

\begin{Prop}
  \begin{enumerate}
  \item
$(a,g)\circ \mR(t)=\mR(t)\circ (a^{-1},g)$ for all $(a,g)\in GL_1\times G_1.$
  \item
 For any $f\in \mS_c$ and $w\in C_0$ there exists $\epsilon$ such that
    for $|t|<\epsilon$ holds
    $$\mR(t)(f)(w)=\mR(0)(f)(w)+ c_{\psi,q}\chi_K(t)|t|^{n-2}\mR_1(f)(w),$$
    where $$\mR_1(f)(w)=\int\limits_F f(sw)\chi_K(-s)|s|^{n-2}ds$$ 
\end{enumerate}
\end{Prop}

\begin{proof}
  \begin{enumerate}

\item
For $(a,g)\in GL_1\times G_1$ holds:
  $$
  \mR(t)((a,g) f)(w)=\chi_K(a)|a|^{n-1}
  \int\limits_{Y_w(t)} f(avg) |\omega_{w,t}(v)|=$$
  $$\chi_K(a)|a|^{1-n}
  \int\limits_{Y_{w}(t)} f(avg) |\omega_{a^{-1}wg,t}(avg)|=
  \chi_K(a)|a|^{1-n}
  \int\limits_{Y_{a^{-1}wg}(t)} f(v) \omega_{a^{-1}wg,t}(v)=$$
 $$ \chi_K(a^{-1})|a|^{1-n}\mR(t)(f)(a^{-1}wg)=(a^{-1},g)\mR(t)(f)(w).$$
  \item
For any $f\in \mS_c$ and $w\in C_0$  there exists $\epsilon$ such that for
$|t|<\epsilon$ holds
$$\mR(t)(f)(w)=\int\limits_{F}\delta_{V_w}(st)(f_{t,s}) ds =
\int\limits_{F}\delta_{V_w}(st)(f_{0,s}) ds, $$
which by \ref{delta:difference} equals
$$\int\limits_F \delta_{V_w}(0)(f_{0,s}) ds + c_{\psi,q}
\chi_K(t)|t|^{n-2}\int\limits_F f_{0,s}(0)\chi_K(s)|s|^{n-2} ds=$$
$$\mR(0)(f)(w)+ c_{\psi,q}\chi_K(t)|t|^{n-2}
\int\limits_{F} f(-sw)\chi_K(s)|s|^{n-2}ds=$$
$$\mR(0)(f)(w)+ c_{\psi,q}\chi_K(t)|t|^{n-2}\mR_1(f)(w)$$
as required. 
\end{enumerate}
  \end{proof}
  
\subsection{The normalization}
The normalized Radon transform $\hat \mR:\mS_c\rightarrow \mS^\infty$
is defined by
\begin{equation}\label{def:norm:radon}
  \hat\mR(f)(w)=\int\limits_F \mR(t)(f)(w)\psi(t) dt.
  \end{equation}
By eqation (\ref{Radon:homog})  it follows that
$$\hat\mR(f)(xw)=\int\limits_F \mR(t)(f)(w)\psi(xt) dt.$$
The following crucial properties of $\hat\mR$ follow
easily from the properties of $\mR(t):$
\begin{Prop}\label{hatR:asymp}
  For any $w\in C_0,$ and $f\in \mS_c$ holds 
  \begin{enumerate}
  \item $\hat \mR\circ (a,g)=(a^{-1},g)\circ \hat \mR$ for all
    $(a,g)\in GL_1    \times G_1.$
 % \chi_K(a)|a|^{1-n} \hat R(f)(a^{-1}w)$.
\item  The germ  $[\hat \mR(f)]_0=I_C(f)$ is constant.
\item The function $x\mapsto \hat \mR(f)(xw) \cdot |x|^{n-1}$ is bounded
  as $|x|\to \infty$.
 
  \end{enumerate}
\end{Prop}

\begin{proof}
  \begin{enumerate}
  \item The first part follows from the equivariance property
    of the operator $\mR(t)$ for any $t$. 
    \item
    For a fixed $f$ and a compact neighborhood $W$ of $0$ in $C$
    the support of the function $x\mapsto \mR(x)(f)(w)$ is uniformly bounded
    for all $w\in W$. Hence for $|t|$ small enough and $w\in W$ holds
    $\mR(f)(x)(w)\psi(tx)=\mR(f)(x)(w)$. This implies
  $$\hat \mR(f)(tw)=\int\limits_F \mR(x)(f)(w)\psi(tx) dx=
  \int\limits_F \mR(x)(f)(w) dx=\int\limits_C f(v)\omega(v)$$
  for $w\in W$. Hence $[\hat\mR(f)]_0=I_C(f)$ as required.

\item 
  For any $f\in \mS_c$ and $w\in C_0$ the function
  $\mR(\cdot)(f)(w)$ is a sum of a smooth  function of compact support
  on $F$ and
  a function $\phi$ such that $\phi(s)=c \cdot\chi_K(s)|s|^{n-2}$ for $|s|\le 1$ and $0$ otherwise. 
A standard computation shows that $\mF_\psi(\phi)(x)|x|^{n-1}$ is bounded. 
  \end{enumerate}
\end{proof}

\subsection{The Fourier transform}

We have all the ingredients to define  the  Fourier transform on the cone,
mentioned in the introduction. It is an operator 
$\Phi:\mS_c\rightarrow \mS^\infty.$

\begin{Def} 
For any  $f\in \mS_c$ define for each $w \in C_0$
$$\Phi(f)( w)=
\gamma(\chi_K,\psi)
\int\limits_{F^\times}\hat \mR(f)(x w)
\psi(x^{-1}) \chi_K(-x) |x|^{n-2} d^\times x.$$
\end{Def}

%In next section we shall show that the
%restriction of the operator $\Pi(r)$ to $\mS_c$ equals to $\Phi$.
%This implies that  the image of $\Phi$ is contained in $\mS$. 

The integral converges absolutely due to properties of the function
$x\mapsto \hat \mR(f)(xw)$, proved in \ref{hatR:asymp}, parts $(2)$ and $(3)$.
The operator $\Phi:\mS_c\rightarrow \mS^\infty$ satisfies
$$(a,g)\circ \Phi=\Phi\circ (a^{-1},g), (a,g)\in GL_1\times G_1,$$
since $\hat \mR$ satisfies the same property.

In  section \ref{sec:Pi(r)} we shall show that the
restriction of the operator $\Pi(r)$ to $\mS_c$ equals to $\Phi$.
This implies that  the image of $\Phi$ is contained in $\mS$. 

There is a natural decomposition
   \begin{equation}\label{Fourier:decomp}
    \Phi=\Phi_1+\Phi_2,\end{equation}
 where the operators
  $\Phi_1, \Phi_2: \mS_c\rightarrow \mS^\infty$ are defined by 
\begin{equation}\label{def:Phi1}
 \Phi_1(f)(w)= \gamma(\chi_K,\psi)\int\limits_{F^\times}
 \hat \mR(f)(xw)\chi_K(-x) |x|^{n-2} d^\times x.
 \end{equation}
%leading term

\begin{equation}\label{def:Phi2}
  \Phi_2(f)(w)=
\gamma(\chi_K,\psi)\int\limits_{F}
\hat \mR(f)(xw) (\psi(x^{-1})-1)\chi_K(-x) |x|^{n-2} d^\times x.
\end{equation}
%constant term

\begin{Prop}\label{images:Phi12}
  \begin{enumerate}
    \item
  The image of $\Phi_1$ equals to $\mS_1\subset \mS(|\cdot|^{-1})$.
\item
For any $f\in S_c$ one has $[\Phi_2(f)]_0=c_{\psi,q} I_C(f).$
In particular,  the image of $\Phi_2$ is contained in the space of
  locally constant functions on $C$.
\end{enumerate}
  \end{Prop}

\begin{proof}
\begin{enumerate}\item  
  It is clear from the formula that $\Phi_1$
  is well-defined and belongs to
$\Hom_{GL_1\times G_1}(\mS_c, |\cdot|^{-1}\otimes \mS(|\cdot|^{-1}))$
that is one-dimensional. 
In particular, by \ref{mult:one:Phi1}
$\Phi_1(f)\in |\cdot|\otimes \Pi_1=\mS_1.$

\item
  The support of the function
  $(\psi(x^{-1})-1)\chi_K(-x) |x|^{n-2}$ is a bounded set in $F$, denote
  it $\mathcal B$. By \ref{hatR:asymp} $(2)$
 there is a neighborhood $W$  of $0$ in $C$ such that for
  $\hat \mR(f)(xw)=I_C(f)$ for $w\in W, x\in \mathcal B$.
%constant term
Hence  $[\Phi_2(f)]_0=c_{\psi,q} I_C(f)$ for all $f\in \mS_c,$ as required.    
\end{enumerate}
\end{proof}

\section{The operator $\Pi(r)$}\label{sec:Pi(r)}
In introduction we have defined an involution $r$ in $O(V)$.
In this section we prove
\begin{Thm}\label{Pi(r)=Phi}
  The restriction of the operator $\Pi(r)$ to $\mS_c$ equals $\Phi.$
\end{Thm}  
The proof  occupies the entire section. 
As the proof is partly computational, it is instructive first
to sketch the idea behind the proof.

Recall the decomposition  $V=\bH\oplus \bH \oplus V_2,$
and  the  bases $\{e,e^\ast\}$,  $\{e_1,e_1^\ast\}$
of the first and the second copy of $\bH$ respectively.

Define two involutions $r_1, r_2\in O(V)$ such that
$$r_1(e)=e_1, \quad r_1(e_1)=e, \quad r_1|_{V_2}=Id$$
$$r_2(e)=e, \quad r_2(e_1)=e_1^\ast, \quad r_2|_{V_2}=Id$$

It is easy to see that $r=r_1\cdot r_2\cdot r_1$.
Since $r_2\in M$, the formula for the  operator $\Pi(r_2)$
is known by \ref{intro:Q:action}.

\begin{itemize}
  \item
In subsection \ref{subsec:Pi(r1)} we derive a formula 
for the action of $\Pi(r_1)$ on a subspace of $\mS_c$, dense in $L^2(C,\omega)$.
To do this we use a mixed model for the minimal representation
realized on a space of functions on $F^\times \times V_2$.
In this model the Heisenberg parabolic subgroup $P=LU$ acts explicitly.
Since $r_1\in L$, it is easy to write a formula for its action on
$\mS_c(F^\times\times V_2)$.

The map
$$ F^\times \times  V_2\hookrightarrow C_0, (y,w)\mapsto (-y,w,q(w)/y)$$
gives rise to an isometric isomorphism  of the Hilbert spaces 
$$L^2(F\oplus  V_2, d^\times x\wedge dv)\simeq L^2(C,\omega).$$
We show that it is $G$-equivariant. In particular, the formula
for $\Pi(r_1)$ can be written on the image of
$\mS_c(F^\times\times V_2)$, that is dense subspace in $L^2(C,\omega_C)$. 

\item
We compute $\Pi(r)=\Pi(r_1)\circ \Pi(r_2)\circ \Pi(r_1)$ on this subspace. 
The result is an integral over $F$ involving 
 values of the distribution $\delta_{V_2}(t)\circ \mF_{\psi, V_2}$.
The decomposition in \ref{delta:difference} gives rise to a decomposition
$\Pi(r)=T_1+T_2$. It remains to show that $T_1=\Phi_1$ and
$T_2=\Phi_2$. 
\item
  We verify directly that $T_2=\Phi_2$ on the domain of $T_2$.
  The equivariance properties of $T_1$ imply that it
  can be extended to an operator on $\mS_c$, proportional to $\Phi_1$.
  Hence $\Pi(r)=c \Phi_1+\Phi_2$ or equivalently
  $\Pi(r)-\Phi=(c-1)\Phi_1$. The left hand side has image
  in the space of functions of bounded support, while the right hand side
  in the space of homogeneous functions, unless $c=1$. Hence $c=1$ and
  $\Pi(r)=\Phi$.
  \end{itemize}

 \subsection{The operator $\Pi(r_1)$} \label{subsec:Pi(r1)}
\subsubsection{The structure of Heisenberg parabolic subgroup}
 Let $P=LU$ be the Heisenberg parabolic subgroup of $G$.
 One has $L\simeq GL_2\times O(V_2)$.

On the space $W=\Span\{e,e_1\}$ we define
a symplectic structure by fixing $[e,e_1]=1$.

The unipotent radical $U$ is isomorphic to the Heisenberg
group associated to the symplectic space $W\otimes V_2$
with the form $(\cdot,\cdot)=[\cdot,\cdot]\otimes q$. 
As a set it is $W\otimes V_2\oplus F$ with the group law
$$(x_1,t_1)\cdot (x_2,t_2)=(x_1+x_2, t_1+t_2+\frac{1}{2}(x_1,x_2)).$$
We denote by $l(x)$ the element $(x,0)$ in $U$.

Let us further give notation for elements in $GL_2=GL(W)$. 
 
$$t_y=\left(\begin{array}{cc} 1&0 \\0 &y \end{array}\right),\quad
t^y=\left(\begin{array}{cc} y&0 \\0 &1 \end{array}\right),\quad
t(y)=\left(\begin{array}{cc} y& 0\\0 &y^{-1} \end{array}\right),\quad
,$$
$$n(b)=\left(\begin{array}{cc} 1&b \\0 &1 \end{array}\right),\quad
w_0=\left(\begin{array}{cc}0 &1 \\-1 & 0\end{array}\right).$$

  \subsubsection{ The Weil representation}
  The derived group $P'$ of $P$ that is
  isomorphic to $SL(W)\times O(V_2)\rtimes U$
  acts via Weil representation $\omega_{\psi,q}$ on the space $\mS_c(V_2).$
  Let us record the action of certain elements. The center $Z$ of $U$ acts by the character $\psi$.

\begin{equation}\label{weil:formula}
  \left\{\begin{array}{ll}
    \omega_{\psi,q}(n(b)) f(v)=\psi(b q(v))f(v)& b\in F\\
    \omega_{\psi,q}(t(y)) f(v)=\chi_K(y)|y|^{n-1}f(yv)& y\in F^\times\\
    \omega_{\psi,q}(w_0) f(v)=
    \gamma(\chi_K,\psi)\int\limits_{V_2}f(u)\psi(\la u,v\ra) du\\
\omega_{\psi,q}(h)f(v)=f(vh)& h\in O(V_2)\\
    \omega_{\psi,q}(l(e\otimes u)) f(v)=\psi(-\la u,v\ra)f(v) & u\in V_2\\
    \omega_{\psi,q}(l(e_1\otimes u)) f(v)=f(v+u) & u\in V_2
    \end{array}
  \right.
 \end{equation}
  
  The space of the compact induction $\ind^P_{P'}\omega_{\psi,q}$
  can be identified with $\mS_c(F^\times\times V_2)$ via
$$\beta:\ind^P_{P'}\mS_c(V_2)\rightarrow \mS_c(F^\times\times V_2),
\quad \beta(f)(y,v)=f(t_y)(v)$$

Using the formula above we can write the action of $r_1$ on
$\mS_c(F^\times\times V_2)$.

\begin{Lem}\label{r_1:mixed}
  The action of the element $r_1$ is given by
  $$r_1 (f)(y,w)=
 % \gamma(\chi_K,\psi)\chi_K(-y)|y|^{n-1}\int\limits_{V_2} \cdot f(-y,-yu)\psi(\la u,v\ra) du$$
\gamma(\chi_K,\psi)\int\limits_{V_2} t^{-y}\cdot f(1,u)\psi(\la u,v\ra) du$$

\end{Lem}

\begin{proof}
For any $p\in P'$  one has by definition
$$pt_b\cdot  f(y,v)=\omega_{\psi,q}(t_ypt_y^{-1})f(by,v).$$

Note that   $r_1=w_0t(-1)t_{-1}$ and $t_y w_0t(-1) t_y^{-1}= w_0t(-y).$

  $$r_1(f)(y,v)=\omega_{\psi,q}(w_0 t(-y))f(-y,v)=
  \gamma(\chi_K,\psi)\int\limits_{V_2} \omega_{\psi,q}(t(-y))f(-y,u)\psi(\la u,v\ra) du=$$ 
  $$\gamma(\chi_K,\psi)\int\limits_{V_2} \chi_K(-y)|y|^{n-1}f(-y,-yu)\psi(\la u,v\ra) du=
\gamma(\chi_K,\psi)\int\limits_{V_2} t^{-y}f(1,u)\psi(\la u,v\ra) du. $$
\end{proof}

\subsubsection{The mixed model}
The action of $P$ preserves the inner product on $\mS_c(F^\times\times V_2)$
with respect to the measure $dx^\times \wedge dv$. 
By unitary completion of the $P$-module $\mS_c(F^\times\times V_2)$
we obtain the unitary action of $P$ on the Hilbert space $L^2(F\oplus V_2).$
By the pioneering work of \cite{KazhdanSavin}, see also \cite{GanSavin}
this action of $P$ can be extended uniquely to a unitary representation of $G$,
that is isomorphic to the unitary completion $\hat \Pi$ of $\Pi$.
We call the space  $L^2(F\oplus V_2)$ the mixed model of
the minimal representation and denote the action by $\hat \Pi^m$,
to distinguish it from $\hat \Pi$, which is the realization
of the minimal representation in Schrodinger model $L^2(C)$.
%The restriction of $\hat \Pi^m$ to the space of smooth vectors is denoted by $\Pi^m$.

It is easy to  write the transition between the models. 

Define an embedding with open dense image
$$\alpha:F^\times\times V_2\hookrightarrow C_0, \quad \alpha(y,v)=(y,v,-q(v)/y)$$
and denote the image by $C_0^1$. 

This gives rise to the isomorphism
$$\alpha^\ast: \mS_c(C_0^1)\rightarrow \mS_c(F^\times\times V_2)$$
that can be extended to the isometry 
$$\hat\alpha^\ast: L^2(C_0)\rightarrow L^2(F\oplus  V_2).$$

\begin{Prop} The map $\hat \alpha^\ast$ is $G$ equivariant. 
\end{Prop}
\begin{proof}  It is straightforward to check using the formulae \ref{weil:formula}
  and \ref{intro:Q:action} that  $\alpha^\ast$ is $P\cap Q$ equivariant.
  Hence $\hat \alpha^\ast$ is $P\cap Q$ equivariant.
  By Mackey theory $\hat \Pi$ restricted to $P\cap Q$ is irreducible and hence
  $\Hom_{P\cap Q}(\hat \Pi, \hat \Pi^m)$ is one-dimensional. 
  The non-zero space  $\Hom_G(\hat \Pi,\hat \Pi^m)$ is contained
 in ${\Hom_{P\cap Q}(\hat \Pi, \hat \Pi^m)},$ and hence is equal to it.
In particular  $\hat \alpha^\ast$ is $G$-equivariant 
\end{proof}

We use $\alpha$ to write the formula for $\Pi(r_1)$ on $\mS_c(C^1_0)$.
Using Lemma \ref{r_1:mixed} we obtain

\begin{Cor}
  The operator $\Pi(r_1)$ preserves $\mS_c(C^1_0)$
  and the action is given by
$$\Pi(r_1)(f)(y,w,-q(w)/y)=\gamma(\chi_K,\psi)
\int\limits_{V_2} (-y)\cdot f(1,u)\psi(\la u,v\ra) du.$$
\end{Cor}

\subsection{The action of $\Pi(r)$}

Our goal is to use the formula for $\Pi(r_1)$ to write the formula
of $\Pi(r)$ on a dense subspace of $\mS_c$. Below we define this subspace. 

\subsubsection{The space $\mS^2\subset \mS$}
Introduce notation
for open dense subsets on the cone and the spaces of functions
of compact support on them. 
$$C_0^1=\{(y,w,-q(w)/y), y\in F^\times, w\in V_2\}\subset C_0, \quad
\mS^1=\mS_c(C^1_0).$$
Similarly
$$C_0^{2n}=\{(-q(w)/y,w,y), y\in F^\times, w\in V_2\}, \quad
\mS^{2n}=\mS_c(C^{2n}_0).$$
Further put  $\mS^2=\mS^1\cap \Pi(r_1)(\mS^{2n}).$

\begin{Lem}
  \begin{enumerate}
    \item
    The operator $\Pi(r)$ preserves $\mS^2$.
  \item The space $\mS^2$ is dense in $L^2(C)$.
  \end{enumerate}
\end{Lem}

\begin{proof}
  \begin{enumerate}
    \item
  The operator $\Pi(r_1)$ preserves $\mS^1$ and hence
  maps $\mS^2$ to $\mS^1\cap \mS^{2n}$ that is preserved by $\Pi(r_2)$. 
So $\Pi(r)=\Pi(r_1)\circ \Pi(r_2)\circ \Pi(r_1)$ preserves the space $\mS^2.$
\item
The space $\mS^1\cap \mS^{2n}=\mS_c(C_0^1\cap C_0^{2n})$ 
is dense in $L^2(C)$. The space $\mS^2$ contains $\Pi(r)(\mS^1\cap \mS^{2n})$
so is dense in $L^2(C)$.
  \end{enumerate}
  \end{proof}

\subsubsection{The formula for $\Pi(r)$ on $\mS^2$}
For a function $f$ on a vector space $U$ and a vector $u_0\in U$
we denote by $f_{u_0}$ the translation of $f$ by $u_0$, i.e. 
$f_{u_0}(u)=f(u+u_0)$.

\begin{Prop}\label{pi(r)}
  Let $f=\sum_i f_i'\otimes f_i''\in \mS^2\subset\mS_c(F^\times)\otimes \mS_c(V_2).$
One has  
  $$\Pi(r)(f)(\tilde w)=\sum_i\int\limits_F f_i'(s y)
\cdot
\delta_{V_2}(-s^{-1})(\mF_{\psi,V_2}(f^{''}_{i,sw})) \chi_K(-s)|s|^{n-2} d^\times s$$
for $\tilde w=(y,w,-\frac{q(w)}{y}) \in C_0^1$ and is zero otherwise.
\end{Prop}

\begin{proof}
  The operator
  $\Pi(r)=\Pi(r_1)\circ \Pi(r_2)\circ \Pi(r_1)$ preserves $\mS^2$.
  We apply formula for $\Pi(r_1)$ on $\mS^1$ and $\Pi(r_2)$ on $\mS_c$.
We shall write $\Pi(h)(f)$ just as $h\cdot f$ to ease notation.

Applying the formula for $\Pi(r_1)$ and $\Pi(r_2)$ we get
$$r\cdot f(\tilde w)= \int\limits_{V_2}\int\limits_{V_2}
t^{q(u)} \cdot t_{-y}\cdot f (1,v,-q(v)) \psi(\la u,v+w\ra)du dv$$
By Fubini theorem for  the variable $u\in V_2$ we obtain
$$ \int\limits_{F}\int\limits_{V_2(s)}
\int\limits_{V_2} (t^{s}t_{-y}\cdot f)
(1,v,-q(v)) \psi(\la u,v+w\ra)dv |\eta_s(u)|  ds=$$ 
 $$ \int\limits_{F}\int\limits_{V_2(s)}
  \int\limits_{V_2} f (-sy,sv,sq(v)/y)|s|^{n-1}\chi_K(s)
  \psi(\la u,v+w\ra )dv |\eta_s(u)| ds =$$ 
By change of variables $sv\mapsto v$ this equals
  $$\int\limits_F  \int\limits_{V_2(s)}
  \int\limits_{V_2} f (-sy,v,q(v)/sy)\psi(\la s^{-1}u,v+sw\ra) dv |\eta_s(u)|
  \chi_K(s) |s|^{2-n}d^\times s.$$  
  The change of variables $s^{-1}u\mapsto u, v+sw\mapsto v, s\mapsto -s$
  produces
$$ \int\limits_F
  \int\limits_{V_2(-s^{-1})}\int\limits_{V_2}
  f (sy,v+sw,-q(v+sw)/sy)\psi(\la u,v\ra)dv|\eta_{s^{-1}}(u)|
  |s|^{n-2}\chi_K(-s) d^\times s =$$
  $$\sum_i\int\limits_F f_i'(sy)
  \delta_{V_2}(-s^{-1})
  (\mF_{\psi,V_2}(f^{''}_{i,sw})) |s|^{n-2}\chi_K(-s)d^\times s $$
  as required.
\end{proof}

\subsection{The decomposition of $\Pi(r)$}
We use the proposition \ref{delta:difference} and the formula
\ref{H:def} to decompose the operator $\Pi(r)$. 

\begin{Prop}\label{cor:T1}
  The restriction of $\Pi(r)$ to $\mS^2$ decomposes as 
  $$\Pi(r)=T_1+\Phi_2,$$
  where $T_1\in\Hom_{GL_1}(\mS^2,|\cdot|^{-1}\otimes \mS(|\cdot|^{-1}))$.
\end{Prop}

\begin{proof}
  We use the decomposition \ref{H:def}
  to write $\Pi(r)$ as $T_1+T_2$ where
\begin{equation}\label{T1:def}
  T_1(f)(\tilde w)=
  \sum_i\int\limits_F  f_i'(ty) \cdot
  \delta_{V_2}(0)(\mF_{\psi,V_2}f^{''}_{i,tw})\cdot 
  \chi_K(-t)|t|^{n-2} d^\times t
  \end{equation}
  and
  $$T_2(f)=
  \sum_i \int\limits_F f_i'(ty)\int\limits_{V_2} f^{''}_{i,tw}(v) H_{-t^{-1}}(v) dv
  \chi_K(-t)(\psi(t)-1) |t|^{n-2} d^\times t.$$

Let us show that $T_2(f)=\Phi_2(f)$ for all $f\in \mS^2$. 
Plugging the formula for $H_{-t^{-1}}(v)$ from \ref{H:def},
the expression $\gamma(\chi_K,\psi)^{-1} T_2(f)(\tilde w)$  equals 
\begin{equation}\label{T2:unfolded}
\int\limits_F  
\left(\int\limits_F\int\limits_{V_2}f(ty,v+tw,-\frac{q(v+tw)}{ty} )
 \psi(-\frac{q(v)x}{t})dv d^\times t \right)
(\psi(x^{-1})-1)\chi_K(-x)|x|^{n-2} d^\times x.
\end{equation}
To show that $T_2=\Phi_2$ it is enough to see that
the inner integral equals to $\hat \mR(f)(xw)$.

\begin{Lem}\label{T2:to:Phi2}
  Let $\tilde w=(y,w,-\frac{q(w)}{y})\in C_0$ and
  $f\in \mS^1$.
  Then one has
  $$\int\limits_{F^\times} \int\limits_{V_2}
  f(ty, v+tw, -\frac{q(v+tw)}{ty}) \psi(-\frac{q(v)x}{t})dv d^\times t=
  \hat \mR(f)(x \tilde w).$$
\end{Lem}

\begin{proof}
 By Fubini theorem the integral becomes 
$$
\int\limits_{s\in F}
\int\limits_{t\in F}
\int\limits_{V_2(s)}f(ty,v+tw,-\frac{q(v+tw)}{ty}) |\eta_s(v)|
 \psi(-\frac{sx}{t})  d^\times t ds .$$
 
 Let $\tilde v=(ty, v+tw,-\frac{q(v+tw)}{ty})$ and
 $\tilde w=(y,w,-\frac{q(w)}{y}).$
Then
$$\la \tilde v, \tilde w\ra=
-q(w)t-\frac{q(v+tw)}{t}+(v+tw,w)=$$
$$-\frac{\la w,w\ra t}{2}-\frac{\la v+tw,v+tw \ra}{2}+\la v,w\ra +t\la w,w\ra=-t^{-1}q(v).$$
In particular, for $v\in V_2(s)$, one has $\la \tilde v, \tilde w\ra= -s/t$.
The integral above becomes 
$$\int\limits_F 
\int\limits_{s\in F} \int\limits_{v\in V_2(s)}f(\tilde v) 
 \psi(\la \tilde v, \tilde w \ra x) |\eta_s(v)| ds d^\times t=
\int\limits_{C} f(\tilde v)  \psi(\la \tilde v, \tilde w\ra x) |\omega(\tilde v)| =$$
$$\int\limits_{s\in F}
\int\limits_{\varphi_{\tilde w}^{-1}(s)}
f(\tilde v)|\eta_{\tilde w,s}(\tilde v)|\psi(sx) ds=\hat \mR(f)(x\tilde w).$$
\end{proof}
\vskip 5pt
The image of the operator $T_1$ is contained in $\mS(|\cdot|^{-1})$. Indeed, 
$$T_1(f)(a\tilde w)= T_1(f)(ay,aw,-\frac{q(aw)}{ay})=$$
$$\sum_i\int\limits_F f_i'(tay)
\delta_{V_2}(0)(\mF_{\psi,V_2}(f^{''}_{i,taw})) |t|^{n-2}\chi_K(-t)d^\times t=
\chi_K(a)|a|^{2-n}T_1(\tilde w).$$

It remains to show that $GL_1$ equivariance properties of $T_1$.
Let $a\in GL_1$ acting on $\mS$ as in \ref{intro:Q:action}.

  By definition $\Pi(r)\circ \Pi(a)=\Pi(a^{-1})\circ \Pi(r)$.
  By the properties of Radon transform $T_2\circ a=a^{-1}\circ T_2$.
  Hence the same is true for $T_1$. 
  $$T_1(a\cdot f)=a^{-1}\cdot T_1(f)=$$
  $$\chi_K(a)|a|^{n-1} \sum_i  f_i'(a^{-1}y)\int\limits_F
\delta_{V_2}(0)(\mF_{\psi,V_2}(f^{''}_{i,ta^{-1}w}))
  |t|^{n-2}\chi_K(-t)d^\times t= 
  |a|^{-1} T_1(f)(\tilde w)$$
as required.
  \end{proof}

 \subsection{The proof of  Theorem \ref{Pi(r)=Phi}}
 \begin{proof}   
   Recall the  decomposition $\Pi(r)=T_1+\Phi_2$ on $\mS^2$ from
   Proposition \ref{cor:T1}.
The operator
$\tilde \Phi_1=\Pi(r)-\Phi_2$  is a $GL_1\times G_1$ equivariant operator
on $\mS_c$ whose restriction to $\mS^2$  equals $T_1$.

The  $GL_1\times G_1$-equivariance of $\tilde \Phi_1$
and Proposition \ref{cor:T1} together imply that   
$$\tilde \Phi_1\in
\Hom_{GL_1\times G_1}(\mS_c, |\cdot|^{-1}\otimes \mS(|\cdot|^{-1})).$$
On the other side   $\Phi_1$ belongs to the same one-dimensional
space. Hence  there exists a constant $c$ such that $\tilde \Phi_1=c\Phi_1$.
Equivalently,   $\Pi(r)-\Phi=(c-1)\Phi_1$. 

For any $f\in \mS_c$ the function $(\Pi(r)-\Phi)(f)$ is of bounded support
and $\Phi_1(f)\in \,\mS(|\cdot|^{-1})$ is not, unless it is zero.
Hence the equality of the operators is possible only for $c=1$.
So $c=1$ and $\Pi(r)=\Phi$ on $\mS_c,$ as required. 
\end{proof}

\bibliographystyle{alpha}
 \bibliography{bib}

\newcommand{\etalchar}[1]{$^{#1}$}
\begin{thebibliography}{DS{\etalchar{+}}99}

\bibitem[BK99]{BravermanKazhdan}
Alexander Braverman and David Kazhdan.
\newblock On the {S}chwartz space of the basic affine space.
\newblock {\em Selecta Math. (N.S.)}, 5(1):1--28, 1999.

\bibitem[DS{\etalchar{+}}99]{HarishChandra}
Stephen DeBacker, Paul~J Sally, et~al.
\newblock {\em Admissible Invariant Distributions on Reductive $ p $-adic
  Groups}.
\newblock Number~16. American Mathematical Soc., 1999.

\bibitem[GK22]{GurevichKazhdanBAS}
Nadya Gurevich and David Kazhdan.
\newblock Fourier transforms on the basic affine space of a quasi-split group.
\newblock 2022.

\bibitem[GS05]{GanSavin}
Wee~Teck Gan and Gordan Savin.
\newblock On minimal representations definitions and properties.
\newblock {\em Representation Theory of the American Mathematical Society},
  9(3):46--93, 2005.

\bibitem[HM79]{HoweMoore}
Roger~E. Howe and Calvin~C. Moore.
\newblock Asymptotic properties of unitary representations.
\newblock {\em J. Functional Analysis}, 32(1):72--96, 1979.

\bibitem[KM11]{KobayshiMano}
Toshiyuki Kobayashi and Gen Mano.
\newblock The {S}chr\"{o}dinger model for the minimal representation of the
  indefinite orthogonal group {${\rm O}(p,q)$}.
\newblock {\em Mem. Amer. Math. Soc.}, 213(1000):vi+132, 2011.

\bibitem[KS90]{KazhdanSavin}
D.~Kazhdan and G.~Savin.
\newblock The smallest representation of simply laced groups.
\newblock {\em Festschrift in honor of I. I. Piatetski-Shapiro on the occasion
  of his sixtieth birthday, Part I (Ramat Aviv, 1989), Israel Math. Conf.
  Proc., Vol. 2, Weizmann, Jerusalem}, page 209–223, 1990.

\bibitem[MW87]{MoeglinWaldspurger}
Colette M{\oe}glin and Jean-Loup Waldspurger.
\newblock Modeles de whittaker d{\'e}g{\'e}n{\'e}r{\'e}s pour des groupes
  p-adiques.
\newblock {\em Mathematische Zeitschrift}, 196(3):427--452, 1987.

\bibitem[Sav94]{Savin}
Gordan Savin.
\newblock Dual pair {$G_{J}\times{\rm PGL}_2$} where {$G_{J}$} is the
  automorphism group of the {J}ordan algebra {${ J}$}.
\newblock {\em Invent. Math.}, 118(1):141--160, 1994.

\bibitem[SW07]{SavinWoodbury}
Gordan Savin and Michael Woodbury.
\newblock Structure of internal modules and a formula for the spherical vector
  of minimal representations.
\newblock {\em Journal of Algebra}, 312(2):755--772, 2007.

\bibitem[Wei64]{Weil}
Andr\'{e} Weil.
\newblock Sur certains groupes d'op\'{e}rateurs unitaires.
\newblock {\em Acta Math.}, 111:143--211, 1964.

\bibitem[Wei03]{Weissman}
Martin Weissman.
\newblock The fourier-jacobi map and small representations.
\newblock {\em Representation Theory of the American Mathematical Society},
  7(13):275--299, 2003.

\end{thebibliography}

\end{document}